\documentclass[12pt]{amsart}

\usepackage{amsmath,amssymb,amsthm}
\usepackage{graphicx}
\usepackage{hyperref}
\usepackage{booktabs}
\usepackage{enumitem}

\newtheorem{theorem}{Theorem}
\newtheorem{conjecture}{Conjecture}

\newtheorem{remark}{Remark}

\title{Maximal Volume Ideal Polyhedra and the Arithmetic Angle Phenomenon}

\author{Igor Rivin}

\date{\today}
\address{Mathematics Department, Temple University}
\email{igor@dimensionreducers.com}
\subjclass{52A55, 11F06, 52B10}
\keywords{hyperbolic geometry, volume, ideal polyhedra, arithmeticity, FUchsian Groups}
\thanks{The author would like to thank Turing Enterprises for their support, Lambda Labs for a compute tinme grant, and Claude Opus 4.5 for coding assistance. The author would also like to thank A.~Kolpakov for interesting conversations.}

\begin{document}

\maketitle

\begin{abstract} 
We present a software suite for the analysis and optimization of ideal convex polyhedra in hyperbolic 3-space $\mathbb{H}^3$. Using Rivin's variational characterization of ideal polyhedra, we develop efficient algorithms for checking combinatorial realizability and finding volume-maximizing configurations. Our systematic computational study reveals two striking phenomena: (1) maximal volume ideal polyhedra consistently exhibit dihedral angles that are rational multiples of $\pi$---a property with no obvious explanation from the optimization formulation; and (2) the distribution of volumes for random configurations is well-approximated by a Beta distribution, with mean normalized volume converging to approximately $\ln 2 \approx 0.69$ as the vertex count increases. We provide complete data for small vertex counts, including vertex positions, triangulations, and verified rational angle structures. An interactive implementation is publicly available.
\end{abstract}

\section{Introduction}

Ideal polyhedra in hyperbolic 3-space $\mathbb{H}^3$ occupy a central position in low-dimensional topology and geometry. Their study connects to the geometry of hyperbolic 3-manifolds, Teichmüller theory, and the representation theory of surface groups. The volume of hyperbolic 3-manifolds, a topological invariant, can often be computed by decomposing the manifold into ideal tetrahedra---the simplest ideal polyhedra, and the variational methods developed by the author form the foundation of the so-called Casson-Rivin program, which is central to modern study of geometry and topology of 3-dimensional manifolds (see, for example, \cite{futer2009angled}.

A natural question arises: among all ideal polyhedra with a given number of vertices, which has the largest volume? This optimization problem, while analytically straightforward (it reduces to maximizing a concave function over a convex polytope), reveals surprising structure. The maximal configurations exhibit a form of ``arithmeticity''---their dihedral angles are consistently rational multiples of $\pi$. Moreover, the distribution of volumes for random configurations follows a Beta distribution, with parameters scaling linearly in the vertex count and mean normalized volume converging to approximately $\ln 2$.

In this paper, we present a comprehensive software suite for:
\begin{enumerate}[label=(\roman*)]
    \item Checking whether a given triangulation is realizable as the boundary of an ideal convex polyhedron
    \item Finding optimal (volume-maximizing) geometric realizations of realizable triangulations
    \item Computing the Fuchsian group associated with the surface of the polyhedron
    \item Visualizing polyhedra in multiple models of hyperbolic space
\end{enumerate}

Our computational experiments reveal a surprising phenomenon: when optimizing for maximal volume, the resulting polyhedra consistently have dihedral angles that are \emph{exact rational multiples of $\pi$}. This arithmetic structure emerges from a purely analytic optimization problem with no built-in number-theoretic constraints.

\section{Background and Notation}

\subsection{Ideal Polyhedra in $\mathbb{H}^3$}

An \emph{ideal polyhedron} in hyperbolic 3-space $\mathbb{H}^3$ is a convex polyhedron whose vertices lie on the boundary at infinity $\partial_\infty \mathbb{H}^3 \cong \mathbb{CP}^1$. Despite having vertices ``at infinity,'' such polyhedra have finite volume---a consequence of the exponential shrinking of hyperbolic space toward the boundary.

The combinatorial structure of an ideal polyhedron is encoded by a triangulation of the 2-sphere: the \emph{Delaunay triangulation} of the ideal vertices as seen from the center of mass of the polyhedron. Conversely, not every triangulation of the sphere arises from an ideal polyhedron---realizability imposes constraints on the possible angle structures, as characterized by Rivin's theorem below.

We work primarily in the upper half-space model, where $\mathbb{H}^3 = \{(x,y,z) \in \mathbb{R}^3 : z > 0\}$ with metric $ds^2 = (dx^2 + dy^2 + dz^2)/z^2$. The boundary $\partial_\infty \mathbb{H}^3$ is the extended complex plane $\mathbb{C} \cup \{\infty\}$, and ideal vertices are specified as points in $\mathbb{CP}^1$.

\subsection{Rivin's Variational Characterization}

The foundational result for our approach is Rivin's characterization of ideal polyhedra:

\begin{theorem}[Rivin \cite{Rivin1994,RivinCombOpt}]
A triangulation $\mathcal{T}$ of the sphere is realizable as the boundary of a convex ideal polyhedron in $\mathbb{H}^3$ if and only if there exist positive real numbers $\theta_c > 0$ assigned to each corner $c$ of each triangle such that:
\begin{enumerate}[label=(\roman*)]
    \item \textbf{Triangle constraint:} For each triangle $T$ with corners $c_1, c_2, c_3$:
    \[
    \theta_{c_1} + \theta_{c_2} + \theta_{c_3} = \pi
    \]
    \item \textbf{Vertex constraint:} For each vertex $v$, summing over all corners at $v$:
    \[
    \sum_{c \text{ at } v} \theta_c = 2\pi
    \]
    \item \textbf{Edge constraint:} For each edge $e$, let $\alpha_e$ and $\beta_e$ be the angles opposite to $e$ in the two triangles sharing $e$. Then:
    \[
    \alpha_e + \beta_e < \pi
    \]
\end{enumerate}
\end{theorem}

The angles $\theta_c$ are the \emph{dihedral angles} of the ideal polyhedron---the angle between adjacent faces along their common edge, measured inside the polyhedron. The key insight is that these constraints are all linear (or strict linear inequalities), so realizability can be checked via linear programming, and volume optimization becomes a concave maximization problem over a convex polytope.

\subsection{The Lobachevsky Function and Hyperbolic Volume}

The hyperbolic volume of an ideal tetrahedron with dihedral angles $\alpha, \beta, \gamma$ (opposite pairs equal) is given by:
\[
V(\alpha, \beta, \gamma) = \Lambda(\alpha) + \Lambda(\beta) + \Lambda(\gamma)
\]
where $\Lambda(\theta)$ is the Lobachevsky function:
\[
\Lambda(\theta) = -\int_0^\theta \log |2\sin t| \, dt
\]

The Lobachevsky function satisfies $\Lambda(\pi - \theta) = -\Lambda(\theta)$ and has maximum value $\Lambda(\pi/3) \approx 0.5074$ (Catalan's constant divided by 2). The function is strictly concave on $(0, \pi)$, which ensures that the volume function over the Rivin polytope has a unique maximum. The regular ideal tetrahedron, with all dihedral angles equal to $\pi/3$, achieves volume $V_4 = 3\Lambda(\pi/3) \approx 1.0149$.

\subsection{The Penner-Rivin Algorithm}

Given a realized ideal polyhedron with its dihedral angles, we can compute the associated Fuchsian group $\Gamma < \mathrm{PSL}(2,\mathbb{R})$ using the Penner-Rivin algorithm \cite{RivinTriangulations}. The surface of the ideal polyhedron inherits a complete hyperbolic metric, making it a hyperbolic surface $S$ with cusps at the ideal vertices. The fundamental group $\pi_1(S)$ acts on the hyperbolic plane via a discrete faithful representation---the \emph{holonomy representation}.

The algorithm proceeds as follows:
\begin{enumerate}[label=(\roman*)]
    \item \textbf{Decorate the triangulation:} Assign to each edge $e$ the complex number
    \[
    z_e = \exp(i(\alpha_e + \beta_e))
    \]
    where $\alpha_e, \beta_e$ are the angles opposite $e$ in the adjacent triangles. This is the \emph{shape parameter} of the edge.

    \item \textbf{Compute edge matrices:} For each oriented edge from vertex $u$ to vertex $v$, construct a $2 \times 2$ matrix encoding the hyperbolic translation along the edge, parameterized by the shape parameters of the surrounding edges.

    \item \textbf{Compose around loops:} For each generator of $\pi_1(S)$, represented as a path crossing edges of the triangulation, multiply the corresponding edge matrices to obtain the holonomy matrix in $\mathrm{PSL}(2,\mathbb{R})$.
\end{enumerate}

The resulting matrices generate the Fuchsian group $\Gamma$. When all shape parameters have modulus 1 (equivalently, when all $\alpha_e + \beta_e$ are real), the group is \emph{arithmetic} if certain additional number-theoretic conditions hold. Our observation that maximal volume configurations have rational angles $p\pi/q$ suggests a connection to arithmetic Fuchsian groups, since in these cases $z_e = \exp(ip\pi/q)$ is a root of unity.

\section{Methods}

\subsection{Realizability Testing}

Given a triangulation $\mathcal{T}$ with $n$ vertices, $f$ faces (triangles), and $e$ edges, we check realizability by solving a linear feasibility problem. Introduce variables $\theta_c$ for each corner $c$ (there are $3f$ corners total). The constraints from Theorem 1 become:

\textbf{Linear program for realizability:}
\begin{align*}
\text{Find } & \theta_c \text{ such that:} \\
\theta_{c_1} + \theta_{c_2} + \theta_{c_3} &= \pi && \text{for each triangle } T \\
\sum_{c \text{ at } v} \theta_c &= 2\pi && \text{for each vertex } v \\
\alpha_e + \beta_e &\leq \pi - \epsilon && \text{for each edge } e \\
\theta_c &\geq \epsilon && \text{for each corner } c
\end{align*}
where $\epsilon > 0$ is a small tolerance (we use $\epsilon = 10^{-6}$) to convert strict inequalities to non-strict ones. If this LP is feasible, the triangulation is realizable.

\subsection{Volume Optimization}

The hyperbolic volume of an ideal polyhedron is the sum of volumes of its constituent ideal tetrahedra:
\[
V(\boldsymbol{\theta}) = \sum_T \left( \Lambda(\theta_{c_1}) + \Lambda(\theta_{c_2}) + \Lambda(\theta_{c_3}) \right)
\]
where the sum is over all triangles $T$ and $\Lambda$ is the Lobachevsky function. Since $\Lambda$ is strictly concave and the constraint set is convex, this is a concave maximization problem with a unique global maximum.

We use a two-phase optimization strategy:
\begin{enumerate}[label=(\roman*)]
    \item \textbf{Global search:} Differential evolution with population size $15 \times d$ (where $d$ is the number of free parameters) explores the feasible region to find promising basins of attraction.
    \item \textbf{Local refinement:} L-BFGS-B polishes the best solution found, exploiting the smoothness of the Lobachevsky function.
\end{enumerate}

For robustness, we run 100 independent trials per vertex count and report the best solution found.

\subsection{Rational Angle Detection}

To detect whether an angle $\theta$ is a rational multiple of $\pi$, we compute the continued fraction expansion of $\theta/\pi$:
\[
\frac{\theta}{\pi} = a_0 + \cfrac{1}{a_1 + \cfrac{1}{a_2 + \cdots}}
\]
Truncating at each level yields convergents $p_k/q_k$ that are best rational approximations. If some convergent satisfies $|\theta/\pi - p_k/q_k| < 10^{-10}$ with $q_k \leq 100$, we report $\theta = p_k\pi/q_k$ as a rational angle.

In all maximal volume configurations tested, this procedure successfully identifies exact rational multiples---the numerical angles match the rational predictions to machine precision.

\section{Results}

\subsection{Maximal Volume Configurations}

We computed maximal volume ideal polyhedra for $n = 4, 5, \ldots, 12$ vertices using 100 independent optimization trials per vertex count. Table~\ref{tab:maximal-volumes} summarizes our findings.

\begin{table}[htbp]
\centering
\caption{Maximal volume ideal polyhedra with $n$ vertices. $V_4 \approx 1.0149$ is the volume of the regular ideal tetrahedron. The ``Faces'' column gives $2n-4$ by Euler's formula.}
\label{tab:maximal-volumes}
\begin{tabular}{ccccc}
\toprule
$n$ & Volume & Faces & $V/V_4$ & Angle denominator $q$ \\
\midrule
4 & 1.014942 & 4 & 1.0000 & 3 \\
5 & 2.029883 & 6 & 2.0000 & 3 \\
6 & 3.663862 & 8 & 3.6099 & 2 \\
7 & 4.986773 & 10 & 4.9134 & 5 \\
8 & 6.488469 & 12 & 6.3929 & 11 \\
9 & 8.162538 & 14 & 8.0424 & mixed \\
10 & 9.839315 & 16 & 9.6946 & 14 \\
11 & 11.449290 & 18 & 11.2810 & 7, 14 \\
12 & 13.529628 & 20 & 13.3308 & 5 \\
\bottomrule
\end{tabular}
\end{table}

\begin{remark}[Comparison with Platonic solids]
For $n=6$, the ``obvious'' candidate is the ideal right-angled octahedron, with all dihedral angles equal to $\pi/2$. Indeed, our optimization confirms this: the maximal configuration has angle denominator $q=2$, meaning all angles are multiples of $\pi/2$. The right-angled octahedron is optimal among 6-vertex ideal polyhedra.

For $n=8$, the natural candidate is the ideal cube (with vertices at the 8 corners of a cube inscribed in the sphere). The ideal cube has $V/V_4 = 5$, giving volume $V_{\text{cube}} \approx 5.07$. However, our optimization finds a configuration with $V \approx 6.49$, exceeding the cube by approximately 28\%. The optimal 8-vertex polyhedron is \emph{not} combinatorially a cube---it has a different triangulation with the surprising angle denominator $q=11$.

For $n=12$, the natural candidate is the ideal regular icosahedron. The icosahedron has $V/V_4 = 10$, giving $V_{\text{icosahedron}} \approx 10.15$. Our optimization finds $V \approx 13.53$, exceeding the icosahedron by approximately 33\%. As with the cube, the optimal 12-vertex polyhedron has a different combinatorial structure. Remarkably, all dihedral angles are multiples of $\pi/5$: specifically $2\pi/5$, $3\pi/5$, $4\pi/5$, and $\pi$ (a flat edge).
\end{remark}

\begin{remark}[Symmetry of optimal configurations]
The automorphism groups of the optimal triangulations show a clear pattern: high symmetry for small $n$, rapidly decreasing thereafter.
\begin{center}
\begin{tabular}{ccc}
$n$ & $|\mathrm{Aut}|$ & Structure \\
\hline
4 & 6 & $S_3$ (tetrahedron) \\
5 & 6 & $D_3$ (bipyramid) \\
6 & 8 & $D_4$ (octahedron) \\
7 & 4 & $\mathbb{Z}_2 \times \mathbb{Z}_2$ \\
8--12 & 2--4 & $\mathbb{Z}_2$ or $\mathbb{Z}_2 \times \mathbb{Z}_2$
\end{tabular}
\end{center}
Volume maximization does not favor symmetric configurations---the optimal polyhedra for $n \geq 7$ have only a single reflection symmetry or none at all.
\end{remark}

\begin{figure}[htbp]
\centering
\includegraphics[width=\textwidth]{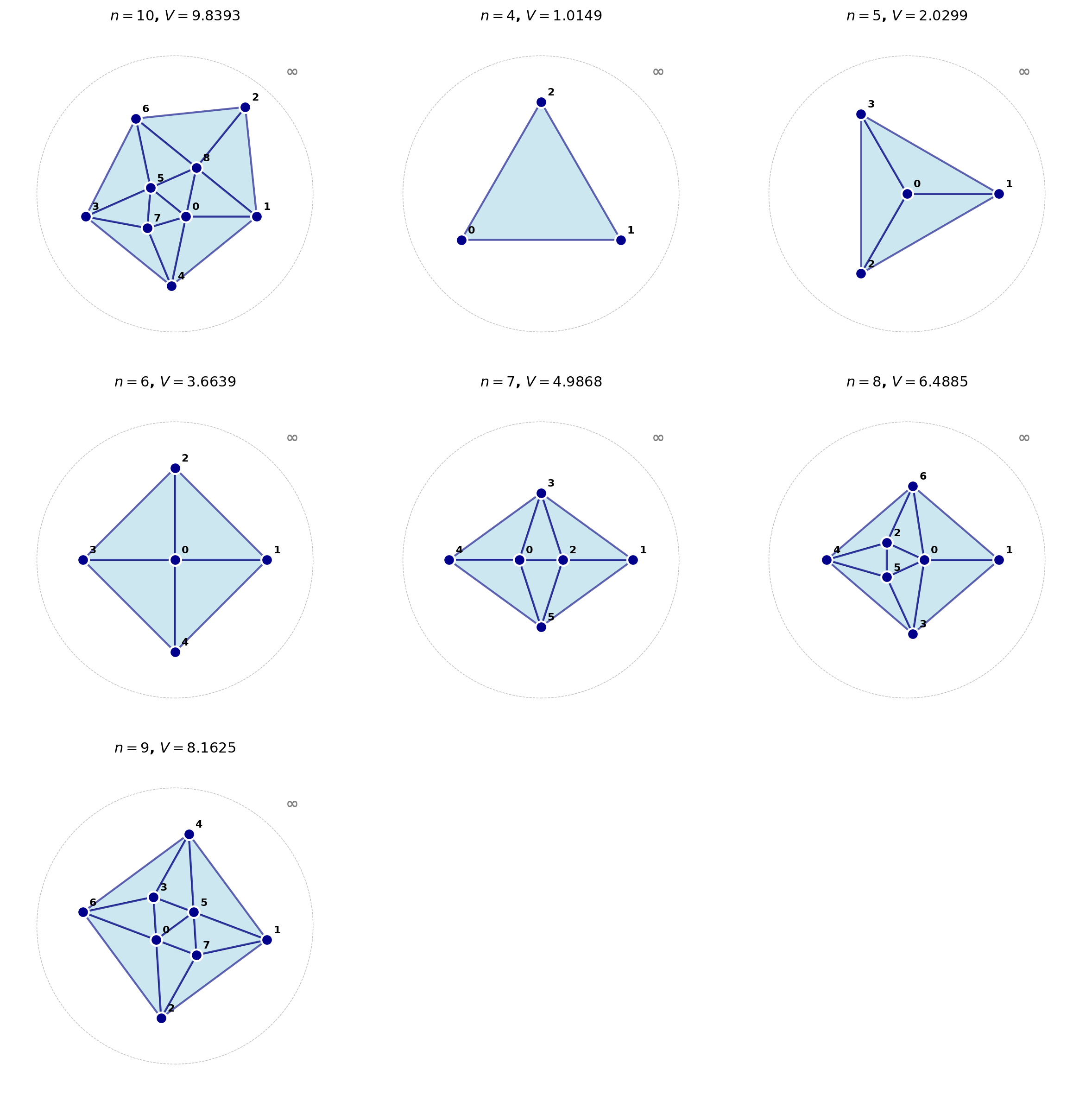}
\caption{Planar Delaunay triangulations of the maximal volume ideal polyhedra. Each triangulation represents the combinatorial structure of the polyhedron projected onto the plane, with the point at infinity as the unbounded face.}
\label{fig:triangulations}
\end{figure}

\begin{figure}[htbp]
\centering
\includegraphics[width=\textwidth]{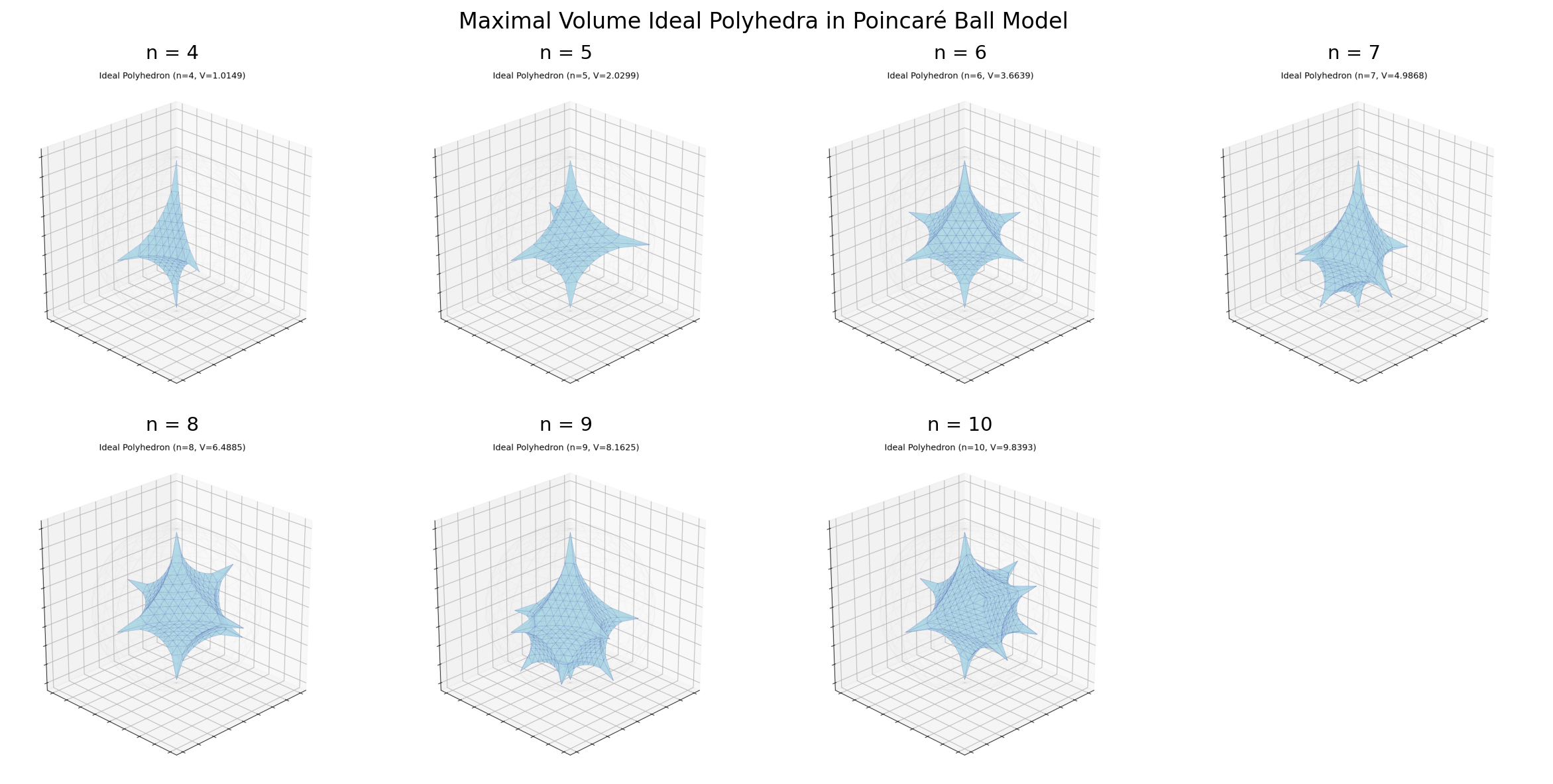}
\caption{Maximal volume ideal polyhedra visualized in the Poincar\'e ball model of $\mathbb{H}^3$. Vertices lie on the boundary sphere $\partial_\infty \mathbb{H}^3$, and the geodesic faces curve inward as characteristic of hyperbolic space.}
\label{fig:poincare}
\end{figure}

\subsection{The Rational Angle Phenomenon}

Our most striking finding is that volume-maximizing configurations \emph{within each combinatorial type} exhibit dihedral angles that are exact rational multiples of $\pi$. This holds not only for the globally maximal configurations shown in Table~\ref{tab:maximal-volumes}, but for every realizable triangulation we have tested.

\begin{conjecture}[Rational Angle Conjecture]
For every realizable triangulation $\mathcal{T}$ of the sphere, the volume-maximizing ideal polyhedron with combinatorial type $\mathcal{T}$ has all dihedral angles of the form $\frac{p\pi}{q}$ for integers $p, q$.
\end{conjecture}

This phenomenon is remarkable because:
\begin{itemize}
    \item The optimization has no built-in preference for rational angles
    \item The Lobachevsky function has no special behavior at rational multiples of $\pi$
    \item The constraint set (Rivin polytope) has no obvious arithmetic structure
    \item The phenomenon persists across \emph{all} combinatorial types, not just the globally optimal ones
\end{itemize}

For the globally maximal configurations shown in Table~\ref{tab:maximal-volumes}, we verified the rational angle property in all cases $n = 4, \ldots, 12$. The denominators range from small values like $q = 3$ (tetrahedron), $q = 2$ (octahedral $n=6$), and $q = 5$ ($n=7$ and $n=12$) to larger values like $q = 11$ for $n=8$ and $q = 14$ for $n = 10$ and $n=11$.

\subsection{Connection to Arithmetic Groups}

The prevalence of rational angles in maximal volume configurations suggests a deep connection to arithmetic hyperbolic geometry. When the dihedral angles are rational multiples of $\pi$, the shape parameters $z_e = \exp(i(\alpha_e + \beta_e))$ are roots of unity. This places strong constraints on the holonomy representation.

A Fuchsian group $\Gamma < \mathrm{PSL}(2,\mathbb{R})$ is \emph{arithmetic} if it is commensurable with the group of units in a quaternion algebra over a totally real number field. Arithmetic Fuchsian groups are characterized by having trace field $\mathbb{Q}(\mathrm{tr}(\gamma) : \gamma \in \Gamma)$ that is a number field, and satisfying additional ramification conditions.

When all shape parameters are roots of unity, the traces of holonomy matrices lie in cyclotomic fields. While this alone does not guarantee arithmeticity, it is highly suggestive. We conjecture that:

\begin{conjecture}[Arithmetic Holonomy]
For maximal volume ideal polyhedra, the holonomy representation $\rho: \pi_1(S) \to \mathrm{PSL}(2,\mathbb{R})$ has image commensurable with an arithmetic Fuchsian group.
\end{conjecture}

For the Platonic cases, the Fuchsian groups are known principal congruence subgroups:
\begin{itemize}
    \item $n=4$ (tetrahedron): $\Gamma(3)$
    \item $n=6$ (octahedron): $\Gamma(4)$
    \item Icosahedron: $\Gamma(5)$ --- but this is \emph{not} maximal for $n=12$
\end{itemize}
Identifying the Fuchsian groups corresponding to the maximal volume configurations for $n \geq 7$ remains an interesting open problem. The angle denominators (e.g., $q=11$ for $n=8$, $q=14$ for $n=10$) suggest these may be congruence subgroups of higher level.

\begin{remark}[Coxeter polytopes and branched covers]
A striking distinction emerges between the Platonic and non-Platonic cases. For $n=4$ and $n=6$, \emph{all} dihedral angles have the form $\pi/m$ for an integer $m$ (specifically $\pi/3$ and $\pi/2$ respectively). This means these ideal polyhedra are \emph{Coxeter polytopes}---fundamental domains for discrete reflection groups in $\mathbb{H}^3$.

For $n \geq 7$, however, the angles are $k\pi/q$ with $k > 1$. For instance, $n=7$ has angles $2\pi/5$, $3\pi/5$, $4\pi/5$, and $n=8$ has angles $4\pi/11$, $5\pi/11$, etc. These are \emph{not} Coxeter fundamental domains. The common denominator $q$ suggests these polyhedra may arise as fundamental domains for finite-index subgroups of reflection groups, or equivalently, that their universal covers are branched covers over some link complement. This connection to orbifold geometry and branched coverings merits further investigation.
\end{remark}

\subsection{Volume Distribution of Random Configurations}

To understand the landscape of ideal polyhedron volumes, we sampled 5000 random configurations for each vertex count $n \in \{5, 6, 7, 8, 10, 12\}$ by placing $n-3$ uniformly random points on the sphere (in addition to two fixed vertices and infinity) and computing the volume via the Delaunay triangulation.

Normalizing each volume by the maximal volume $V_{\max}(n)$ from Table~\ref{tab:maximal-volumes}, we find that the normalized volume distribution $V/V_{\max}$ is remarkably well-approximated by a Beta distribution:
\[
\frac{V}{V_{\max}} \sim \mathrm{Beta}(\alpha_n, \beta_n)
\]

Table~\ref{tab:beta-fits} summarizes the fitted parameters. The goodness of fit (KS test p-values $> 0.3$ for $n \geq 8$) confirms that the Beta distribution is an appropriate model.

\begin{table}[htbp]
\centering
\caption{Beta distribution fit parameters for normalized volume $V/V_{\max}$.}
\label{tab:beta-fits}
\begin{tabular}{ccccccc}
\toprule
$n$ & $\alpha$ & $\beta$ & Mean & Std & KS stat & p-value \\
\midrule
5  & 2.88  & 1.53  & 0.659 & 0.208 & 0.031 & 0.0001 \\
6  & 6.74  & 4.11  & 0.623 & 0.140 & 0.046 & 0.0000 \\
7  & 9.67  & 4.84  & 0.667 & 0.120 & 0.024 & 0.0053 \\
8  & 13.26 & 6.12  & 0.685 & 0.103 & 0.011 & 0.62 \\
10 & 23.02 & 10.05 & 0.696 & 0.079 & 0.012 & 0.48 \\
12 & 32.56 & 14.49 & 0.692 & 0.067 & 0.014 & 0.32 \\
\bottomrule
\end{tabular}
\end{table}

The parameters exhibit clear linear scaling with $n$:
\[
\alpha_n \approx 4.25n - 19.3, \qquad \beta_n \approx 1.78n - 7.4
\]
with the ratio $\alpha_n/\beta_n \approx 2.39$ remaining approximately constant. This implies that the mean normalized volume
\[
\mathbb{E}\left[\frac{V}{V_{\max}}\right] = \frac{\alpha_n}{\alpha_n + \beta_n} \approx 0.69 \approx \ln 2
\]
stabilizes as $n$ increases. The intriguing appearance of $\ln 2$ (or perhaps a related constant) as the asymptotic mean suggests a deeper connection to hyperbolic volume theory that we have not yet identified.

\begin{figure}[htbp]
\centering
\includegraphics[width=\textwidth]{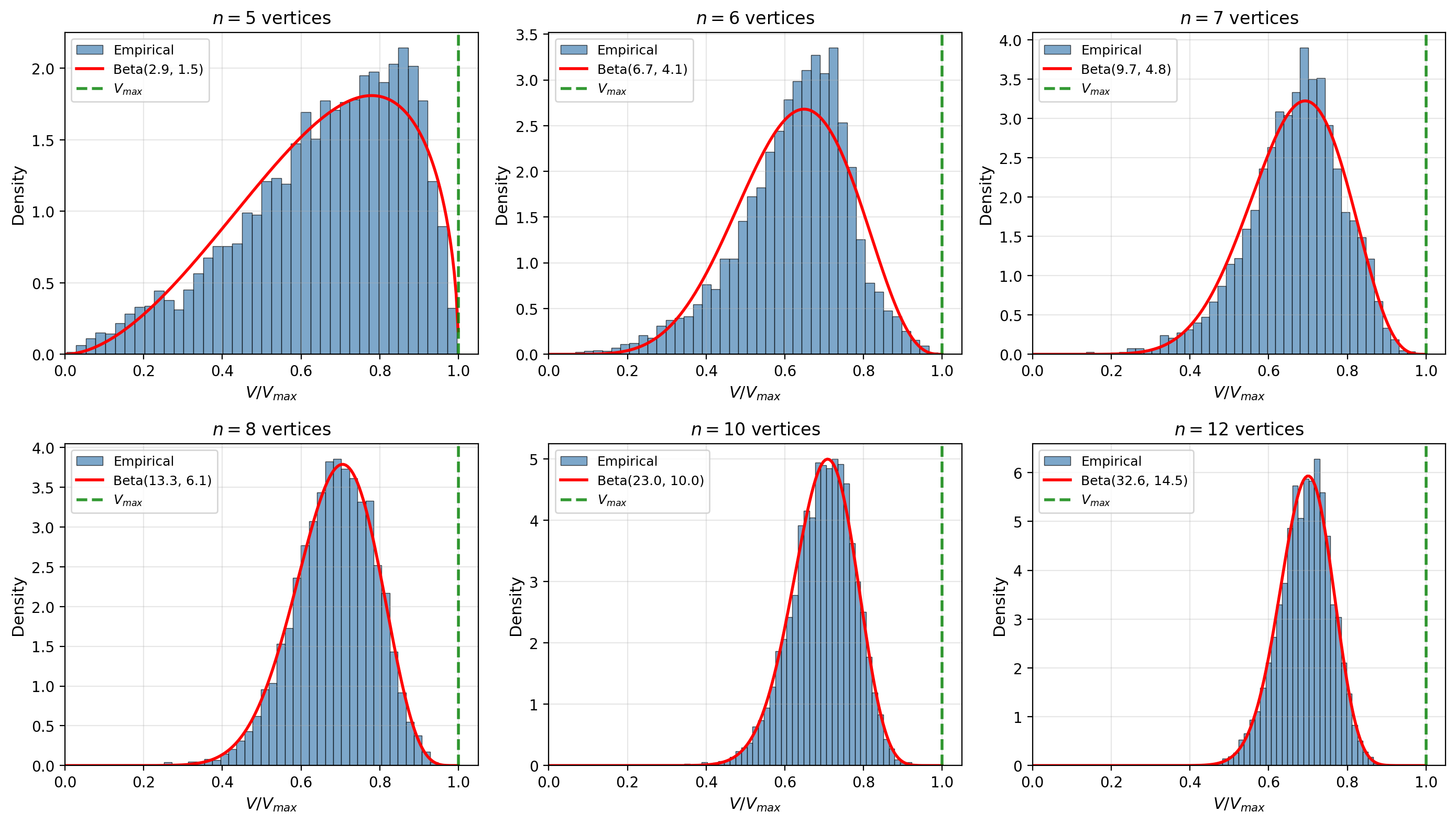}
\caption{Volume distributions of random ideal polyhedra, normalized by $V_{\max}(n)$. Histograms (blue) show empirical distributions from 5000 samples; red curves show fitted Beta($\alpha$, $\beta$) distributions. The fit quality improves with $n$.}
\label{fig:volume-distributions}
\end{figure}

\begin{figure}[htbp]
\centering
\includegraphics[width=\textwidth]{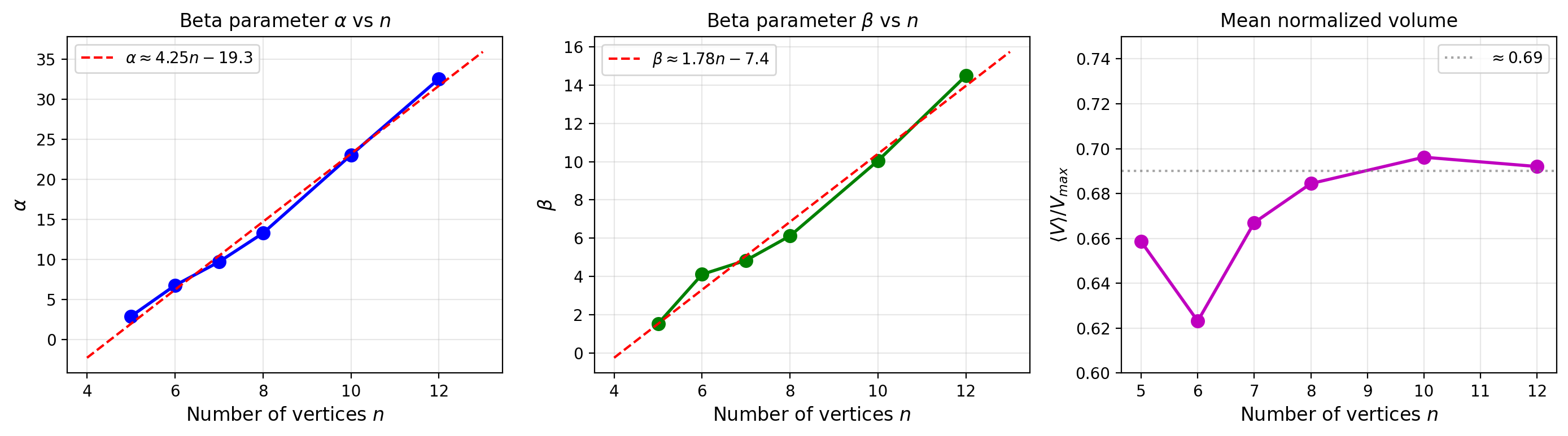}
\caption{Scaling of Beta distribution parameters with vertex count $n$. Left: $\alpha$ grows approximately as $4.25n$. Center: $\beta$ grows as $1.78n$. Right: Mean normalized volume stabilizes near $0.69 \approx \ln 2$.}
\label{fig:beta-scaling}
\end{figure}

The concentration of the distribution (decreasing standard deviation as $n$ increases) indicates that for large $n$, most random configurations achieve volumes close to $69\%$ of the maximum---suggesting that volume optimization yields only modest improvements over random placement for large polyhedra.

\section{Visualizations}

Figure~\ref{fig:triangulations} shows the planar Delaunay triangulations for the maximal volume configurations with $n = 4, \ldots, 12$ vertices. These represent the combinatorial structure of the polyhedra, projected onto the plane via stereographic projection with one vertex (at infinity) as the unbounded outer face.

Figure~\ref{fig:poincare} shows the corresponding ideal polyhedra in the Poincar\'e ball model of $\mathbb{H}^3$. In this model, the ideal vertices lie on the boundary sphere, and the faces appear as geodesic surfaces that curve inward---a characteristic feature of negative curvature. The progression from the simple ideal tetrahedron ($n=4$) to more complex structures is clearly visible.

\section{Software}

The complete software suite is available at:
\begin{itemize}
    \item GitHub: \url{https://github.com/igeal_oly_volume_toolkit}
    \item Interactive demo: \url{https://huggingface.co/spaces/igriv/idealpolyhedra}
\end{itemize}

The implementation includes:
\begin{itemize}
    \item Realizability checking via Rivin's LP characterization
    \item Volume optimization using differential evolution
    \item Fuchsian group computation via the Penner-Rivin algorithm
    \item 3D visualization in Klein, Poincaré ball, and sphere models
    \item Symmetry analysis using nauty
\end{itemize}

\paragraph{Development methodology.} Mathematical software has long suffered from ``bit rot''---the gradual degradation of code as dependencies evolve and environments change. This software was developed using an AI-assisted ``self-maintaining'' approach with Claude Code \cite{RivinDjinn}, which enables automated testing, dependency management, and iterative refinement. This methodology addresses the endemic problem of unmaintained research code by building in mechanisms for autonomous repair and adaptation.

\section{Open Questions}

\begin{enumerate}
    \item \textbf{Proof of Rational Angles}: Can we prove that maximal volume implies rational dihedral angles?

    \item \textbf{Denominator Patterns}: What determines the denominator $q$ for a given $n$? Our data suggests connections to $n-2$ or $n-3$ in some cases.

    \item \textbf{Uniqueness}: Are maximal volume configurations unique up to symmetry?

    \item \textbf{Asymptotics}: What is the growth rate of $V_{\max}(n)$ as $n \to \infty$?

    \item \textbf{Beta Distribution}: Why does the normalized volume follow a Beta distribution? Is the limiting mean exactly $\ln 2$, or merely close to it? What geometric principle underlies the linear scaling of $\alpha_n$ and $\beta_n$ with $n$?
\end{enumerate}



\begin{thebibliography}{99}

\bibitem{futer2009angled}
{Futer, David and Gu{\'e}ritaud, Fran{\c{c}}ois}, ``{Angled decompositions of arborescent link complements}'',
  \emph{Proceedings of the London Mathematical Society},
  \textbf{98}(2) (2009),
325--364


\bibitem{Rivin1994}
I.~Rivin, ``Euclidean structures on simplicial surfaces and hyperbolic volume,'' \emph{Ann. of Math.} \textbf{139} (1994), 553--580.

\bibitem{RivinCombOpt}
I.~Rivin, ``Combinatorial optimization in geometry,'' \emph{Adv. in Appl. Math.} \textbf{31} (2003), 242--271.

\bibitem{RivinTriangulations}
I.~Rivin, ``Triangulations into groups,'' arXiv:math/0510613, 2005.

\bibitem{RivinDjinn}
I.~Rivin, ``I Built a `Djinn' That Fixes My Code While I Sleep,'' \emph{Impure Thoughts} (Substack), December 2025. \url{https://impurethoughts.substack.com/}

\bibitem{Milnor}
J.~Milnor, ``Hyperbolic geometry: The first 150 years,'' \emph{Bull. Amer. Math. Soc.} \textbf{6} (1982), 9--24.

\end{thebibliography}
\end{document}